\documentclass[preprint,12pt,square,comma,sort&compress]{elsarticle}

\usepackage{indentfirst}
\usepackage{mathrsfs}
\usepackage{dsfont}
\usepackage{amssymb}
\usepackage{stmaryrd}
\usepackage{cite}
\usepackage{amsfonts}
\usepackage{graphicx}
\usepackage{multirow}
\usepackage{amsmath}
\usepackage{pifont}
\usepackage{color}
\usepackage{amsthm}
\usepackage{hyperref}

\allowdisplaybreaks[4]

\def\sgn{\mathop{\rm sgn}}
\def\ds{\displaystyle}

\newtheorem{lemma}{Lemma}[section]
\newtheorem{remark}{Remark}[section]

\newtheorem{theorem}{Theorem}[section]

\journal{Applied Mathematics and Computation}

\begin{document}

\begin{frontmatter}

\title{Optimal Harvesting of a Stochastic Lotka-Volterra Competition Model with Periodic Coefficients \tnoteref{myr}}
\tnotetext[myr]{The work was supported by the National Natural Science Foundation of China (No. 12071292, 42450269).}

\author{Wenmin Deng}
\author{Fu Zhang\corref{cor1}}
\address{College of Science, University of Shanghai for Science and Technology, Shanghai 200093, China}

\cortext[cor1]{Corresponding author.}
\ead{fuzhang82@gmail.com}

\begin{abstract}
This paper systematically investigates the optimal harvesting of a stochastic Lotka-Volterra competition model with periodic coefficients. Sufficient conditions for the extinction and persistence in the time average of each species are established. Using Khasminskii's stability theory with suitable Lyapunov functions, we establish sufficient conditions to guarantee the existence of positive periodic solutions to the model. Under certain assumptions, the stability in distribution of this model is proved. Then, we obtain the existence of an optimal harvesting policy and provide explicit expressions for the optimal harvesting effort and the maximum sustainable yield. Finally, we demonstrate our key findings numerically using the Euler-Maruyama method implemented in Python.
\end{abstract}


\end{frontmatter}

\rm
\section{Introduction}
\label{1}

With the rapid development of society and economy, human overexploitation of biological resources (such as overfishing and deforestation) has led to a series of severe ecological problems, including the degradation of ecosystem functions, significant decline in biodiversity, and even species extinction. Against this background,  achieving the sustainable utilization of renewable biological resources has become one of the core issues in the study of population dynamics.

However, real ecosystems are inevitably subject to various random factors.
Climate change, extreme natural disasters, and the unpredictability of human activities may all cause the optimal harvesting strategies derived from deterministic models to fail.
Therefore, employing stochastic models to investigate population dynamics and harvesting strategies has extremely important theoretical and practical significance.

In the research process of optimal harvesting strategies for autonomous stochastic ecosystems (i.e., systems with parameters that do not change with time), scholars have achieved a series of significant breakthroughs. Beddington and May \citep{MR01} were pioneers in examining the optimal harvesting problem of single-species stochastic models, laying a solid theoretical foundation for this field. Alvarez and Shepp \citep{MR02} constructed sufficient conditions for the existence of optimal harvesting strategies in stochastic ecosystems with the analysis of the Hamilton-Jacobi-Bellman equation, providing a key tool for theoretical proof. In addition, Liu and Bai \citep{MR03} obtained the sufficient and necessary conditions for the existence of optimal harvesting strategies by solving the corresponding Fokker-Planck equation for stochastic predator-prey models. However, it should be noted that the method of solving the corresponding Fokker-Planck equation is often challenging to apply to optimal harvesting problems of population models in more complex systems. Subsequently, Liu \citep{MR04} further optimized the analysis method of optimal harvesting strategies for stochastic population models based on the ergodicity of steady-state distributions, which has effectively promoted the expansion of research on more complex population models. For example, Qiu and Deng \citep{MR051} used the ergodic method to investigate the optimal harvesting problem for stochastic competitive populations with discrete time delays and L\^{e}vy jumps. Subsequently, Qiu et al. \citep{MR05} extended this approach to systems with S-type distributed time delays and L\^{e}vy jumps, thereby refining the relevant theoretical framework. These works have jointly constructed the theoretical framework for the study of optimal harvesting strategies in stochastic ecosystems, providing a solid foundation for subsequent research.

Nevertheless, the population dynamics in nature often exhibit more complex spatiotemporal characteristics. Due to the periodic changes of environmental factors such as seasons, the growth rate of populations and the intensity of random disturbances usually have significant time-varying characteristics. Specifically, in the resource-rich spring and summer seasons, populations often show higher growth potential; while in the resource-scarce autumn and winter seasons, populations face greater survival pressures. This periodicity is mainly reflected in the following three aspects: seasonal fluctuations in food resources; the spatiotemporal heterogeneity of habitats and shelters; and periodic changes in environmental temperature. The interplay of these factors makes traditional autonomous models difficult to accurately describe the actual population dynamics.

In existing studies, Liu \citep{MR06} used the stochastic periodic solution of the model as a bridge to conduct the first exploration of the optimal harvesting problem of a stochastic Gompertz model with periodic coefficients (single-species case). However, species in natural ecosystems often form complex interaction networks through competition, symbiosis, and other mechanisms, which makes multi-species models more able to reflect real ecological relationships. Among them, competition models, as basic models describing resource competition among species, occupy an important position in the study of population dynamics. However, due to the theoretical complexity of non-autonomous stochastic systems, there are still obvious gaps in the research on optimal harvesting of competitive systems with periodic coefficients.

This paper aims to establish the theoretical framework for optimal harvesting of stochastic competitive systems with periodic parameters. It mainly promotes existing research from the following two dimensions: Firstly, based on Khasminskii's asymptotic stability theory, by constructing a new Lyapunov function, the existence of positive periodic solutions of the system is rigorously demonstrated; Secondly, combined with ergodic theory, the mathematical characterization of optimal harvesting strategies is established. The research results not only fill the research gap in the optimal harvesting theory of periodic stochastic competition systems but also provide more accurate decision-making basis for actual ecosystem management.

In summary, a stochastic competition model with periodic coefficients can be expressed as follows:
\begin{equation}\label{fxx}
\left\{\begin{array}{rcl}
\ds \mathrm{d}x_{1}(t)&=&x_{1}(t)\Big[r_{1}(t)-h_1-c_{11}x_{1}(t)-c_{12}x_{2}(t)\Big]\mathrm{d}t  \\
&\quad&\ds+\alpha_{1}(t)x_{1}(t)\mathrm{d}B_{1}(t),\\
\ds \mathrm{d}x_{2}(t)&=&x_{2}(t)\Big[r_{2}(t)-h_2-c_{21}x_{1}(t)-c_{22}x_{2}(t)\Big]\mathrm{d}t\\
&\quad&\ds+\alpha_{2}(t)x_{2}(t)\mathrm{d}B_{2}(t),\\
\end{array}
\right.
\end{equation}
where $x_{1}(t)$ and $x_{2}(t)$ stand for the population size of two species, respectively. $r_{i}(t)>0$ is the growth rate of $x_{i}(t)$, $i=1,2$. $h_{i}>0$ represents the harvesting effort of $x_{i}(t)$, $i=1,2$. $c_{ii}>0$ is the intraspecific competition coefficients of $x_{i}(t)$, $i$=1,2; $c_{ij}$ ($i\neq j$; $i,j$=1,2) denotes the interspecific competition rate, respectively. $\alpha_i^2(t)$, $i=1,2$ denote the intensity of the white noise. ${B_i(t)}$ $i=1,2$ are standard independent Brownian motions defined on a complete probability space $(\Omega, \mathcal{F},\mathscr{P})$. The coefficients $r_{i}(t)$, $\alpha_i(t)$, $i,j=1,2$ are continuous $T$-periodic functions.

To simplify notation, we assume $T=1$ throughout this paper. Following Fan and Wang \citep{MR07}, we define the management objective as the annual reward, i.e., maximizing the expected annual-sustainable yield:
$$Y(H)\overset{\triangle}{=}\liminf\limits_{t\rightarrow+\infty}\int_{t}^{t+1}\sum^2_{i=1}\mathbb{E}(h_ix_i(s))ds$$
where $H=(h_1,h_2)^T$ is the harvesting effort. Our goal is to find the optimal harvesting effort (OHE) $H^*$ that maximizes $Y(H)$ under the constraint of species persistence.
\medskip

\rm
\section{Extinction and persistence}
\label{2}

For clarity in subsequent analysis, we introduce the following notation:
\begin{equation*}
\ds x(t)\overset{\triangle}{=}(x_1(t), x_2(t))^T\in R_+^2;
\end{equation*}
\begin{equation*}
\ds b_{i}(t)\overset{\triangle}{=}r_{i}(t)-h_i-\frac{\alpha_{i}^{2}(t)}{2},\quad i=1,2;
\end{equation*}
\begin{equation*}
\ds\Delta_{1}\overset{\triangle}{=}c_{22}\int_{0}^{1}b_{1}(t)\mathrm{d}t -c_{12}\int_{0}^{1}b_{2}(t)\mathrm{d}t,\quad \Delta_{2}\overset{\triangle}{=}c_{11}\int_{0}^{1}b_{2}(t)\mathrm{d}t-c_{21}\int_{0}^{1}b_{1}(t)\mathrm{d}t;
\end{equation*}
\begin{equation*}
\langle f\rangle^{*}\overset{\triangle}{=}\limsup_{t\rightarrow+\infty}\frac{1}{t}\int_{0}^{t}f(s)\mathrm{d}s,\quad
\langle f\rangle_{*}\overset{\triangle}{=}\liminf_{t\rightarrow+\infty}\frac{1}{t}\int_{0}^{t}f(s)\mathrm{d}s;
\end{equation*}
\begin{equation*}
\ds f^{u}\overset{\triangle}{=}\sup_{t\in[0,+\infty)}f(t).
\end{equation*}
\medskip

\begin{lemma}
Model (\ref{fxx}) has a unique global positive solution $x(t)$ almost surely (a.s.). In particular,
\begin{equation}\label{lnt}
\ds\limsup_{t\rightarrow+\infty}\frac{\ln x_i(t)}{\ln t}\leqslant1, \quad a.s., \quad i=1,2.
\end{equation}
\end{lemma}
\medskip

\begin{proof}[$\mathbf{Proof\ of\ Lemma\ 2.1}$] \rm The proof is a special case of Theorem 2.1 and Lemma 3.4 in Li and Mao \citep{MR08} and hence omitted.
\end{proof}
\medskip

\rm
\begin{lemma}
(Liu et al.\citep{MR09} and Xia et al.\citep{MR10}).
Let $z(t)\in C[\Omega\times[0,+\infty),R_+]$, and $\lim_{t\rightarrow+\infty}F(t)/t=0$ $a.s.$
\medskip

\noindent$\rm(\mathfrak{A})$ If there exist some constants $T>0$, $\lambda_0>0$, and $\lambda$ such that for all $t\geqslant T$,
\begin{equation*}
\ln z(t)\leqslant \lambda t-\lambda_0\int_0^tz(s)\mathrm{d}s+F(t),\ a.s.,
\end{equation*}
then
\begin{equation*}
\left\{\begin{array}{l}
\ds\langle z\rangle^*=\limsup_{t\rightarrow+\infty}\frac{1}{t}\int_0^tz(s)\mathrm{d}s\leqslant \lambda/\lambda_0,\quad a.s., \quad if\ \lambda\geqslant0,\\
\ds\lim_{t\rightarrow+\infty}z(t)=0,\quad a.s., \quad if\ \lambda<0.
\end{array}
\right.
\end{equation*}
$\rm(\mathfrak{B})$ If there exist some constants $T>0$, $\lambda_0>0$, and $\lambda>0$ such that for all $t\geqslant T$,
\begin{equation*}
\ln z(t)\geqslant\lambda t-\lambda_0\int_0^tz(s)\mathrm{d}s+F(t),\ a.s.,
\end{equation*}
then
\begin{equation*}
\langle z\rangle_*=\liminf_{t\rightarrow+\infty}\frac{1}{t}\int_0^tz(s)\mathrm{d}s\geqslant \lambda/\lambda_0,\quad a.s.\\
\quad
\end{equation*}
\end{lemma}
\medskip

Before we state our results, we make a assumption in the following:
\medskip

\noindent $\mathbf{Assumption\ 1.}$ Define $\Delta\overset{\triangle}{=}c_{11}c_{22}-c_{21}c_{12}$ and assume that $\Delta>0$.
\medskip

\begin{lemma}
Suppose that Assumption 1 holds, for model (\ref{fxx}),
\medskip

\noindent$\rm(\uppercase\expandafter{\romannumeral1})$ if $\int_{0}^{1}b_{1}(t)\mathrm{d}t<0$ and $\int_{0}^{1}b_{2}(t)\mathrm{d}t<0$, then both $x_1$ and $x_2$ tend to extinction a.s., i.e., $\lim\limits_{t\rightarrow+\infty}x_i(t)=0$, a.s., i=1,2;

\noindent$\rm(\uppercase\expandafter{\romannumeral2})$ If $\int_{0}^{1}b_{1}(t)\mathrm{d}t>0$ and $\int_{0}^{1}b_{2}(t)\mathrm{d}t<0$, then $x_2$ tends to extinction a.s., and
\begin{equation*}
\lim_{t\rightarrow +\infty}\frac{1}{t}\int_{0}^{t}x_{1}(s)\mathrm{d}s=\frac{\int_{0}^{1}b_{1}(t)\mathrm{d}t}{c_{11}}, \quad a.s.
\end{equation*}
\noindent$\rm(\uppercase\expandafter{\romannumeral3})$ If $\int_{0}^{1}b_{1}(t)\mathrm{d}t<0$ and $\int_{0}^{1}b_{2}(t)\mathrm{d}t>0$, then $x_1$ tends to extinction a.s., and
\begin{equation*}
\lim_{t\rightarrow +\infty}\frac{1}{t}\int_{0}^{t}x_{2}(s)\mathrm{d}s=\frac{\int_{0}^{1}b_{2}(t)\mathrm{d}t}{c_{22}}, \quad a.s.
\end{equation*}
\noindent$\rm(\uppercase\expandafter{\romannumeral4})$ If $\int_{0}^{1}b_{1}(t)\mathrm{d}t>0$ and $\int_{0}^{1}b_{2}(t)\mathrm{d}t>0$.

$\rm(\romannumeral1)$ If $\Delta_1>0$ and $\Delta_2<0$, then $x_2$ tends to extinction a.s., and
\begin{equation*}
\lim_{t\rightarrow +\infty}\frac{1}{t}\int_{0}^{t}x_{1}(s)\mathrm{d}s=\frac{\int_{0}^{1}b_{1}(t)\mathrm{d}t}{c_{11}}, \quad a.s.
\end{equation*}

$\rm(\romannumeral2)$ If $\Delta_1<0$ and $\Delta_2>0$, then $x_1$ tends to extinction a.s., and
\begin{equation*}
\lim_{t\rightarrow +\infty}\frac{1}{t}\int_{0}^{t}x_{2}(s)\mathrm{d}s=\frac{\int_{0}^{1}b_{2}(t)\mathrm{d}t}{c_{22}}, \quad a.s.
\end{equation*}

$\rm(\romannumeral3)$ If $\Delta_1>0$ and $\Delta_2>0$, then both $x_1$ and $x_2$ are stable in time average a.s.:
\begin{equation}\label{slu}
\lim_{t\rightarrow +\infty}\frac{1}{t}\int_{0}^{t}x_{1}(s)\mathrm{d}s=\frac{\Delta_{1}}{\Delta}, \quad \lim_{t\rightarrow +\infty}\frac{1}{t}\int_{0}^{t}x_{2}(s)\mathrm{d}s=\frac{\Delta_{2}}{\Delta},\quad a.s.
\end{equation}
\end{lemma}
\medskip

\begin{remark}
Under the conditions $\int_{0}^{1}b_{1}(t)\mathrm{d}t>0$, $\int_{0}^{1}b_{2}(t)\mathrm{d}t>0$ and $\Delta>0$, the inequalities $\Delta_1<0$ and $\Delta_2<0$ cannot be satisfied simultaneously.
\end{remark}
\medskip

\begin{proof}[$\mathbf{Proof\ of\ remark\ 2.1}$] \rm Due to $\Delta_1<0$ and $\Delta_2<0$, we can see
\begin{equation*}
\begin{array}{rcl}
\ds c_{22}\int_{0}^{1}b_{1}(t)\mathrm{d}t< c_{12}\int_{0}^{1}b_{2}(t)\mathrm{d}t,\quad
\ds c_{11}\int_{0}^{1}b_{2}(t)\mathrm{d}t< c_{21}\int_{0}^{1}b_{1}(t)\mathrm{d}t.
\end{array}
\end{equation*}
Note that $c_{11},c_{12},c_{21},c_{22}>0$, $\int_{0}^{1}b_{1}(t)\mathrm{d}t>0$ and $\int_{0}^{1}b_{2}(t)\mathrm{d}t>0$, multiplying both sides of the above two inequalities, we can deduce that
\begin{equation*}
\ds c_{11}c_{22}\int_{0}^{1}b_{1}(t)\mathrm{d}t\int_{0}^{1}b_{2}(t)\mathrm{d}t< c_{12}c_{21}\int_{0}^{1}b_{1}(t)\mathrm{d}t\int_{0}^{1}b_{2}(t)\mathrm{d}t,
\end{equation*}
it follows that
\begin{equation*}
\ds c_{11}c_{22}- c_{12}c_{21}=\Delta<0.
\end{equation*}
The contradiction arises.
\end{proof}
\medskip

\begin{proof}[$\mathbf{Proof\ of\ Lemma\ 2.3}$] \rm Applying the It$\rm\hat{o}$'s formula to model (\ref{fxx}) yields
\begin{equation*}
\begin{array}{rcl}
&\quad&\ds\ln x_{1}(t)-\ln x_{1}(0)\\
&=&\ds \int_{0}^{t}(r_{1}(s)-h_1-\frac{\alpha_{1}^{2}(s)}{2})\mathrm{d}s-c_{11}\int_{0}^{t}x_{1}(s)ds-c_{12}\int_{0}^{t}x_{2}(s)\mathrm{d}s\\
&\quad&\ds +\int_{0}^{t}\alpha_{1}(s)\mathrm{d}B_{1}(t)\\
&=&\ds
\int_{0}^{t}b_{1}(s)\mathrm{d}s-c_{11}\int_{0}^{t}x_{1}(s)\mathrm{d}s-c_{12}\int_{0}^{t}x_{2}(s)\mathrm{d}s+\int_{0}^{t}\alpha_{1}(s)\mathrm{d}B_{1}(s).\\
\end{array}
\end{equation*}
That is to say
\begin{equation}\label{fx1}
\begin{array}{rcl}
\ds \frac{\ln x_{1}(t)}{t}
&=&\ds \frac{\ln x_{1}(0)}{t}+\frac{1}{t}\int_{0}^{t}b_{1}(s)\mathrm{d}s-\frac{c_{11}}{t}\int_{0}^{t}x_{1}(s)\mathrm{d}s\\
&\quad&\ds -\frac{c_{12}}{t}\int_{0}^{t}x_{2}(s)\mathrm{d}s+ \frac{1}{t}\int_{0}^{t}\alpha_{1}(s)\mathrm{d}B_{1}(s).\\
\end{array}
\end{equation}
Similarly, we can deduce that
\begin{equation}\label{fx2}
\begin{array}{rcl}
\ds \frac{\ln x_{2}(t)}{t}
&=&\ds \frac{\ln x_{2}(0)}{t}+\frac{1}{t}\int_{0}^{t}b_{2}(s)\mathrm{d}s-\frac{c_{21}}{t}\int_{0}^{t}x_{1}(s)\mathrm{d}s\\
&\quad&\ds -\frac{c_{22}}{t}\int_{0}^{t}x_{2}(s)\mathrm{d}s+ \frac{1}{t}\int_{0}^{t}\alpha_{2}(s)\mathrm{d}B_{2}(s).\\
\end{array}
\end{equation}
Consider the integral \(\int_0^t b_{i}(s) \, ds\), \(i=1,2\). We can decompose the interval \([0, t]\) into whole numbers of periods and a remainder part. Let \( n \) be the largest integer less than or equal to \( t \), then \( t = n + r \), where \( 0 \leq r < 1 \).
Therefore,
\[
\int_0^t b_{i}(s) \, ds = \int_0^n b_{i}(t) \, dt + \int_n^{n+r} b_{i}(t) \, dt, \quad i=1,2.
\]
Since \( b_{i}(t) \), \(i=1,2\) is a continuous 1-periodic function, we have:
\[
\int_0^n b_{i}(t) \, dt = n \int_0^1 b_{i}(t) \, dt, \quad i=1,2.
\]
Combining the above results, we get:
\[
\int_0^t b_{i}(s) \, ds = n \int_0^1 b_{i}(t) \, dt + \int_n^{n+r} b_{i}(t) \, dt, \quad i=1,2.
\]
Dividing by \( t \) and take limits:
\[
\lim_{t \to +\infty} \frac{1}{t} \int_0^t b_{i}(s) \, ds = \lim_{t \to +\infty} \frac{n}{t} \int_0^1 b_{i}(t) \, dt + \lim_{t \to +\infty} \frac{1}{t} \int_n^{n+r} b_{i}(t) \, dt, \quad i=1,2.
\]
For the remaining part \(\int_n^{n+r} b_{i}(t) \, dt\), since \( r < 1 \), we have:
\[
\left| \int_n^{n+r} b_{i}(t) \, ds \right| \leq \max_{t \in [0, 1]} |b_{i}(t)| \cdot r, \quad i=1,2,
\]
so
\[
 \lim_{t \to +\infty} \frac{1}{t} \int_n^{n+r} b_{i}(t) \, dt=0, \quad i=1,2.
\]
Due to \( n = t - r \), where \( 0 \leq r < 1 \), so \(n/t \to 1\) as \( t \to +\infty \). Therefore:
\[
\lim_{t \to +\infty} \frac{1}{t} \int_0^t b_{i}(s) \, ds = \int_0^1 b_{i}(t) \, dt, \quad i=1,2.
\]
Thus for arbitrary $\varepsilon>0$, there is a random time $T_0>0$ such that for $t\geqslant T_0$,
\begin{equation}\label{zz1}
\ds \int_{0}^{1}b_{i}(t)\mathrm{d}t-\varepsilon\leqslant \frac{1}{t}\int_{0}^{t}b_{i}(s)\mathrm{d}s\leqslant\int_{0}^{1}b_{i}(t)\mathrm{d}t+\varepsilon, \quad i=1,2.
\end{equation}
Note that
\begin{equation*}
\ds \lim_{t\rightarrow+\infty}\frac{\ln x_{i}(0)}{t}=0, \quad i=1,2.
\end{equation*}
Therefore for arbitrary $\varepsilon>0$, there is a random time $T_1>0$ such that for $t\geqslant T_1$,
\begin{equation}\label{zz3}
\ds -\varepsilon\leqslant \frac{\ln x_{i}(0)}{t}\leqslant\varepsilon.
\end{equation}
Furthermore, by the strong law of large numbers (see, e.g., Theorem 3.4 in Chapter 1 of \citep{MR12}), we obtain

\begin{equation}\label{zz2}
\begin{array}{rcl}
\ds \lim_{t\rightarrow+\infty} \frac{1}{t}\int_{0}^{t}\alpha_{i}(s)\mathrm{d}B_{i}(t)=0,\quad a.s., \quad i=1,2.\\
\end{array}
\end{equation}
The derivation of Equation (\ref{zz2}) enables the application of Lemma 2.2 in the following proof.

\noindent $\mathbf{First}$, we prove $\rm(\uppercase\expandafter{\romannumeral1})$. Based on (\ref{fx1}) (\ref{zz1}) and (\ref{zz3}), we get
\begin{equation*}
\frac{\ln x_{1}(t)}{t}
\leqslant \int_{0}^{1}b_{1}(t)\mathrm{d}t+2\varepsilon-\frac{c_{11}}{t}\int_{0}^{t}x_{1}(s)\mathrm{d}s+ \frac{1}{t}\int_{0}^{t}\alpha_{1}(s)\mathrm{d}B_{1}(s).
\end{equation*}
Since $\int_{0}^{1}b_{1}(t)\mathrm{d}t<0$ and let $\varepsilon$ be sufficiently small such that $\int_{0}^{1}b_{1}(t)\mathrm{d}t+2\varepsilon<0$, substituting $(\rm\mathfrak{A})$ of Lemma 2.2, we have
\begin{equation*}
\lim_{t\rightarrow+\infty}x_{1}(t)=0, \quad a.s.
\end{equation*}
In the same way, by (\ref{fx2}), we can demonstrate that if $\int_{0}^{1}b_{2}(t)\mathrm{d}t<0$, then $\lim\limits_{t\rightarrow+\infty}x_{2}(t)=0$, $a.s.$
\medskip

\noindent $\mathbf{Second}$, we prove $\rm(\uppercase\expandafter{\romannumeral2})$. Because $\int_{0}^{1}b_{2}(t)\mathrm{d}t<0$, it follows from $\rm(\uppercase\expandafter{\romannumeral1})$ that $\lim\limits_{t\rightarrow+\infty}x_{2}(t)=0$, $a.s.$ Hence, for every $\varepsilon>0$, there is a random time $T_2>0$ such that for $t\geqslant T_2$,
\begin{equation*}
\ds -\varepsilon\leqslant \frac{c_{12}}{t}\int_{0}^{t}x_{2}(s)\mathrm{d}s\leqslant\varepsilon.
\end{equation*}
By applying the above inequality and equations (\ref{zz1}), (\ref{zz3}) into equation (\ref{fx1}), let $T=\max(T_0,T_1,T_2)$. Then for $t\geqslant T$, one can obtain that
\begin{equation}\label{a}
\begin{array}{rcl}
\ds \frac{\ln x_{1}(t)}{t}
&\leqslant&\ds \int_{0}^{1}b_{1}(t)\mathrm{d}t+3\varepsilon-\frac{c_{11}}{t}\int_{0}^{t}x_{1}(s)\mathrm{d}s+ \frac{1}{t}\int_{0}^{t}\alpha_{1}(s)\mathrm{d}B_{1}(s),\\
\end{array}
\end{equation}
\begin{equation}\label{b}
\begin{array}{rcl}
\ds \frac{\ln x_{1}(t)}{t}
&\geqslant\ds \int_{0}^{1}b_{1}(t)\mathrm{d}t-3\varepsilon-\frac{c_{11}}{t}\int_{0}^{t}x_{1}(s)\mathrm{d}s+ \frac{1}{t}\int_{0}^{t}\alpha_{1}(s)\mathrm{d}B_{1}(s).\\
\end{array}
\end{equation}
Thanks to $\int_{0}^{1}b_{1}(t)\mathrm{d}t>0$ and $\varepsilon$ is arbitrary, we can choose $\varepsilon$ sufficiently small such that $\int_{0}^{1}b_{1}(t)\mathrm{d}t-3\varepsilon>0$. Applying $(\rm\mathfrak{A})$ and $(\rm\mathfrak{B})$ in Lemma 2.2 to (\ref{a}) and (\ref{b}) respectively, we have
\begin{equation*}
\frac{\int_{0}^{1}b_{1}(t)\mathrm{d}t-3\varepsilon}{c_{11}}\leqslant\langle x_{1}\rangle_{*}\leqslant\langle x_{1}\rangle^{*}\leqslant\frac{\int_{0}^{1}b_{1}(t)\mathrm{d}t+3\varepsilon}{c_{11}}.
\end{equation*}
Letting $\varepsilon\rightarrow0$, we obtain that $\lim\limits_{t\rightarrow+\infty}t^{-1}\int_{0}^{t}x_{1}(s)\mathrm{d}s=\int_{0}^{1}b_{1}(t)\mathrm{d}t/c_{11}$, $a.s$.
\medskip

\noindent $\mathbf{Third}$, we prove $(\rm\uppercase\expandafter{\romannumeral3})$. The proof of $(\rm\uppercase\expandafter{\romannumeral3})$ follows by symmetry in a manner analogous to $(\rm\uppercase\expandafter{\romannumeral2})$ and is therefore omitted.
\medskip

\noindent $\mathbf{Fourth}$, we prove $(\rm\uppercase\expandafter{\romannumeral4})$. We consider the following equation:
\begin{equation*}
\left\{\begin{array}{rcl}
\ds \mathrm{d}y_{1}(t)&=&y_{1}(t)\Big[r_{1}(t)-h_1-c_{11}y_{1}(t)\Big]\mathrm{d}t+\alpha_{1}(t)y_{1}(t)\mathrm{d}B_{1}(t),\\
\ds \mathrm{d}y_{2}(t)&=&y_{2}(t)\Big[r_{2}(t)-h_2-c_{22}y_{2}(t)\Big]\mathrm{d}t+\alpha_{2}(t)y_{2}(t)\mathrm{d}B_{2}(t),\\
\end{array}
\right.
\end{equation*}
where $y_{i}(0)=x_{i}(0)$, $i=1,2$. According to the stochastic comparision theorem in Huang \citep{MR11}, we obtain
$$
\ds x_{1}(t)\leqslant y_{1}(t),\quad x_{2}(t)\leqslant y_{2}(t).
$$
Due to $\int_{0}^{1}b_{1}(t)\mathrm{d}t>0$ and $\int_{0}^{1}b_{2}(t)\mathrm{d}t>0$, an argument similar to $(\rm\uppercase\expandafter{\romannumeral2})$ yields
$$\lim_{t\rightarrow+\infty}\frac{1}{t}\int_{0}^{t}y_{i}(s)\mathrm{d}s=\frac{\int_{0}^{1}b_{i}(t)\mathrm{d}t}{c_{ii}},\quad a.s., \quad i=1,2.$$
On the other hand, by computing (\ref{fx2})$\times c_{11}$-(\ref{fx1})$\times c_{21}$, we obtain
\begin{equation}\label{c}
\begin{array}{ll}
\ds\frac{c_{11}}{t}\ln \frac{x_{2}(t)}{x_{2}(0)}
&=\ds \frac{c_{21}}{t}\ln \frac{x_{1}(t)}{x_{1}(0)}+\frac{c_{11}}{t}\int_{0}^{t}b_{2}(s)\mathrm{d}s-\frac{c_{21}}{t}\int_{0}^{t}b_{1}(s)\mathrm{d}s-\frac{\Delta}{t} \int_{0}^{t}x_{2}(s)\mathrm{d}s\\
&\quad\ds +\frac{c_{11}}{t}\int_{0}^{t}\alpha_{2}(s)\mathrm{d}B_{2}(s)-\frac{c_{21}}{t}\int_{0}^{t}\alpha_{1}(s)\mathrm{d}B_{1}(s).\\
\end{array}
\end{equation}
From (\ref{lnt}), for any $\varepsilon>0$,  there exists a random time $T_3>0$ such that for all $t\geqslant T_3$,
\begin{equation}\label{zjb}
\ds\frac{c_{ji}}{t}\ln \frac{x_{i}(t)}{x_{i}(0)}<\varepsilon, \quad \rm for \it  \quad i,j=1,2,\quad i\neq j.
\end{equation}
Substituting (\ref{zz1}), (\ref{zz3}) and (\ref{zjb}) into (\ref{c}), we derive
\begin{equation}\label{d}
\begin{array}{rcl}
\ds\frac{c_{11}}{t}\ln x_{2}(t)
&\leqslant& \ds \Delta_{2}+4\varepsilon-\frac{\Delta}{t}\int_{0}^{t}x_{2}(s)\mathrm{d}s\\
&\quad&\ds +\frac{c_{11}}{t}\int_{0}^{t}\alpha_{2}(s)\mathrm{d}B_{2}(s)-\frac{c_{21}}{t}\int_{0}^{t}\alpha_{1}(s)\mathrm{d}B_{1}(s).\\
\end{array}
\end{equation}
Similarly, by computing (\ref{fx1})$\times c_{22}$-(\ref{fx2})$\times c_{12}$ and applying (\ref{zz1}), (\ref{zz3}) and (\ref{zjb}), we obtain
\begin{equation}\label{e}
\begin{array}{rcl}
\ds\frac{c_{22}}{t}\ln x_{1}(t)
&\leqslant&\ds \Delta_{1}+4\varepsilon- \frac{\Delta}{t}\int_{0}^{t}x_{1}(s)\mathrm{d}s\\
&\quad&\ds +\frac{c_{22}}{t}\int_{0}^{t}\alpha_{1}(s)\mathrm{d}B_{1}(s)-\frac{c_{12}}{t}\int_{0}^{t}\alpha_{2}(s)\mathrm{d}B_{2}(s).\\
\end{array}
\end{equation}
for $t>T'$ $(T'=\max(T_0,T_1,T_3))$ and any $\varepsilon>0$.
\medskip

$\mathbf{(\romannumeral1)}$: Since $\Delta_{2}<0$, we may choose $\varepsilon$ sufficiently small such that $\Delta_{2}+4\varepsilon<0$. Applying $(\rm\mathfrak{A})$ in Lemma 2.2 to (\ref{d}), one conclude that

$$\lim\limits_{t\rightarrow +\infty}x_{2}(t)=0, \quad a.s.$$
The proof of
$$\lim\limits_{t\rightarrow +\infty}t^{-1}\int_{0}^{t}x_{1}(s)\mathrm{d}s=\frac{\int_{0}^{1}b_{1}(t)\mathrm{d}t}{c_{11}},\quad a.s.$$
follows analogously to $(\rm\uppercase\expandafter{\romannumeral2})$ and hence is omitted.
\medskip

$\mathbf{(\romannumeral2)}$: The argument for $\mathbf{(\romannumeral2)}$ is symmetric to that of $\mathbf{(\romannumeral1)}$ and is omitted for brevity.
\medskip

$\mathbf{(\romannumeral3)}$: Given $\Delta_{2}>0$, (\ref{d}) and $(\rm\mathfrak{A})$ in Lemma 2.2 imply
$$\langle x_{2}\rangle^*\leqslant \frac{\Delta_{2}+4\varepsilon}{\Delta},\quad  a.s.$$
By the arbitrariness of $\varepsilon$, let $\varepsilon\rightarrow0$ yields
\begin{equation}\label{f}
\ds\langle x_{2}\rangle^*\leqslant \frac{\Delta_{2}}{\Delta},\quad  a.s.
\end{equation}
Similarly, from (\ref{e}) and $(\rm\mathfrak{A})$ in Lemma 2.2, we derive
\begin{equation}\label{g}
\ds\langle x_{1}\rangle^*\leqslant \frac{\Delta_{1}}{\Delta},\quad  a.s.
\end{equation}

\noindent Let $\varepsilon$ be sufficiently small such that $c_{11}\frac{\Delta_1}{\Delta}-2\varepsilon>0$. Substituting (\ref{zz1}), (\ref{zz3}) and (\ref{f}) into (\ref{fx1}), we obtain for sufficiently large $t$:
\begin{equation*}
\begin{array}{rcl}
\ds \frac{\ln x_{1}(t)}{t}
&\geqslant&\ds \int_{0}^{1}b_{1}(t)\mathrm{d}t-2\varepsilon-\frac{c_{11}}{t}\int_{0}^{t}x_{1}(s)\mathrm{d}s-c_{12}\langle x_{2}(t)\rangle^*+ \frac{1}{t}\int_{0}^{t}\alpha_{1}(s)\mathrm{d}B_{1}(s)\\
&\geqslant&\ds \int_{0}^{1}b_{1}(t)\mathrm{d}t-2\varepsilon-\frac{c_{11}}{t}\int_{0}^{t}x_{1}(s)\mathrm{d}s-c_{12}\frac{\Delta_{2}}{\Delta}+\frac{1}{t}\int_{0}^{t}\alpha_{1}(s)\mathrm{d}B_{1}(s)\\
&=&\ds c_{11}\frac{\Delta_{1}}{\Delta}-2\varepsilon-\frac{c_{11}}{t}\int_{0}^{t}x_{1}(s)\mathrm{d}s+ \frac{1}{t}\int_{0}^{t}\alpha_{1}(s)\mathrm{d}B_{1}(s).\\
\end{array}
\end{equation*}
By $(\rm\mathfrak{B})$ in Lemma 2.2 and the arbitrariness of $\varepsilon$, we conclude:
\begin{equation}\label{h}
\ds\langle x_{1}\rangle_*\geqslant \frac{\Delta_{1}}{\Delta},\quad  a.s.
\end{equation}
Similarly, applying (\ref{zz1}), (\ref{zz3}) and (\ref{g}) to (\ref{fx2}) yields:
\begin{equation*}
\langle x_{2}\rangle_*\geqslant \frac{\Delta_{2}}{\Delta},\quad a.s.
 \end{equation*}
Combining this result with (\ref{f})-(\ref{h}), we establish the almost sure convergence $$\lim\limits_{t\rightarrow+\infty}t^{-1}\int_{0}^{t}x_{1}(s)\mathrm{d}s= \frac{\Delta_{1}}{\Delta},\quad \mathrm{and} \quad \lim\limits_{t\rightarrow+\infty}t^{-1}\int_{0}^{t}x_{2}(s)\mathrm{d}s= \frac{\Delta_{2}}{\Delta}, \quad a.s.$$
\end{proof}
\medskip

\rm
\section{Existence of a positive periodic solution}
\label{3}

In this section, we will obtain the sufficient conditions for the existenve of a nontrixial positive periodic solution of model (\ref{fxx}).
\medskip

Considering a $d$-dimensional stochastic differential equation
\begin{equation}\label{xkj}
dx(t)=f(x(t), t)dt+g(x(t), t)dB(t) \quad \text{on } t \ge t_{0}
\end{equation}
with initial value $x(t_{0}) = x_{0} \in \mathbb{R}^{d}$. Define the differential operator $\mathcal{L}$ associated with (\ref{xkj}) by
$$
\mathcal{L} \overset{\triangle}{=} \frac{\partial}{\partial t} + \sum_{k=1}^{d} f_{k}(x, t) \frac{\partial}{\partial x_{k}} + \frac{1}{2} \sum_{k, j=1}^{d} [g^{T}(x, t) g(x, t)]_{|kj} \frac{\partial^2}{\partial x_{k} \partial x_{j}}.
$$

Before we state our result, we make some assumptions in the following:
\medskip

\noindent$\mathbf{Assumption\ 2.}$

\noindent$(H1)$ \sl $\int_{0}^{1}b_{i}(t)\mathrm{d}t>0$, $i=1,2$.

\noindent$(H2)$ \sl $\Phi_m\overset{\triangle}{=}\sum_{i=1}^{2}\int_{0}^{1}b_{i}(t)\mathrm{d}t-\frac{(r_{m}^{u}-h_m+c_{mm}+c_{nm})^{2}}{4c_{mm}}>2$, $ m,n=1,2, m\neq n$.
\medskip

\begin{lemma}
 If Assumption 2 holds, model (\ref{fxx}) has a $1-$periodic solution.
\end{lemma}
\medskip
\begin{proof}[$\mathbf{Proof\ of\ Lemma\ 3.1}$] \rm We will prove this conclusion using Theorem 3.8 in \citep{MR13}.
Define
\begin{equation*}
V(t, x_{1}, x_{2})\overset{\triangle}{=}\sum_{i=1}^{2}(x_{i}(t)-\ln x_{i}(t)+\omega_{i}(t)),
\end{equation*}
where $\omega^{'}_{i}(t)=b_{i}(t)-\int_{0}^{1}b_{i}(t)dt, i=1,2$.
Then $\omega_{i}(t)$ is a $1-$periodic function. In fact, due to $b_{i}(t)$ is a $1-$periodic function,
\begin{equation*}
\begin{array}{rcl}
\ds \omega_{i}(t+1)-\omega_{i}(t)&=&\ds\int_{t}^{t+1}\omega^{'}_{i}(s)\mathrm{d}s\\
&=&\ds\int_{t}^{t+1}b_{i}(s)\mathrm{d}s-\int_{t}^{t+1}\int_{0}^{1}b_{i}(t)\mathrm{d}t\mathrm{d}s\\
&=&\ds\int_{0}^{1}b_{i}(t)\mathrm{d}t-\int_{0}^{1}b_{i}(t)\mathrm{d}t\\
&=&\ds0
\end{array}
\end{equation*}
It is obvious that
\begin{equation}\label{lv1}
\liminf_{k\rightarrow\infty,(x_{1}, x_{2})\in R_+^2\setminus U_{k}}V(t, x_{1}, x_{2})=\infty, \quad i=1,2.
\end{equation}
where $U_{k}=\{(x_{1}, x_{2}): (x_{1}, x_{2})\in(\frac{1}{k},k)\times(\frac{1}{k},k)\}$.

Applying It$\rm\hat{o}$'s formula, we get
\begin{equation*}
\begin{array}{rcl}
\mathcal{L}V&=&\ds\sum_{i=1}^{2}\big[x_{i}(t)\big(r_{i}(t)-h_i-c_{ii}x_{i}(t)-\sum_{j=1,j\neq i}^{2}c_{ij}x_{j}(t)\big)-\big(b_{i}(t)-c_{ii}x_{i}(t)\\
&\quad&\ds-\sum_{j=1,j\neq i}^{2}c_{ij}x_{j}(t)\big)+\omega_{i}^{'}(t)\big]\\
&\leqslant&\ds\sum_{i=1}^{2}\big[x_{i}(t)\big(r_{i}(t)-h_i-c_{ii}x_{i}(t)\big)-\big(b_{i}(t)-c_{ii}x_{i}(t)-\sum_{j=1,j\neq i}^{2}c_{ij}x_{j}(t)\big)\\
&\quad&\ds+b_{i}(t)-\int_{0}^{1}b_{i}(t)\mathrm{d}t\big]\\
&=&\ds-\sum_{i=1}^{2}c_{ii}x_{i}^{2}(t)+\sum_{i=1}^{2}(r_{i}(t)-h_i+c_{ii})x_{i}(t)+\sum_{i,j=1,j\neq i}^{2}c_{ij}x_{j}(t)\\
&\quad&\ds-\sum_{i=1}^{2}\int_{0}^{1}b_{i}(t)\mathrm{d}t\\
&\leqslant&\ds-\sum_{i=1}^{2}c_{ii}x_{i}^{2}(t)+\sum_{i=1}^{2}(r_{i}^{u}-h_i+c_{ii})x_{i}(t)+\sum_{i,j=1,j\neq i}^{2}c_{ij}x_{j}(t)\\
&\quad&\ds-\sum_{i=1}^{2}\int_{0}^{1}b_{i}(t)\mathrm{d}t.
\end{array}
\end{equation*}
Let $U_{\varepsilon}=\{(x_{1}, x_{2}): (x_{1}, x_{2})\in[\varepsilon,\frac{1}{\varepsilon}]\times[\varepsilon,\frac{1}{\varepsilon}]\}$ is a compact set, where choose $\varepsilon$ small enough to satisfy the following conditions:

\begin{equation*}
\begin{array}{@{}l@{\quad}l@{}} 
(1) &\ds (r_{n}^{u}-h_n+c_{nn})\varepsilon \leqslant \frac{1}{2}\biggl\{\sum_{i=1}^{2}\int_{0}^{1}b_{i}(t)\mathrm{d}t-\frac{(r_{m}^{u}-h_m+c_{mm}+c_{nm})^{2}}{4c_{mm}}\biggr\}, \\
(2) &\ds -\frac{c_{nn}}{2\varepsilon}+K \leqslant -1,
\end{array}
\end{equation*}
where $m,n=1,2$, $m\neq n$ and $K$ is defined in the rest of the proof.

Case 1: For any fixed $n (n=1,2)$, if $0< x_{n}<\varepsilon$, we have
\begin{equation*}
\begin{array}{rcl}
\mathcal{L}V&\leqslant&\ds-\sum_{i=1}^{2}c_{ii}x_{i}^{2}(t)+\sum_{i=1}^{2}(r_{i}^{u}-h_i+c_{ii})x_{i}(t)+\sum_{i,j=1,j\neq i}^{2}c_{ij}x_{j}(t)\\
&\quad&\ds-\sum_{i=1}^{2}\int_{0}^{1}b_{i}(t)\mathrm{d}t\\
&\leqslant&\ds-c_{mm}x_{m}^{2}(t)+(r_{m}^{u}-h_m+c_{mm}+c_{nm})x_{m}(t)+(r_{n}^{u}-h_n+c_{nn})\varepsilon\\
&\quad&\ds-\sum_{i=1}^{2}\int_{0}^{1}b_{i}(t)\mathrm{d}t\\
&=&\ds-c_{mm}\big( x_{m}(t)+\frac{r_{m}^{u}-h_m+c_{mm}+c_{nm}}{2c_{mm}}\big)^{2}+\frac{(r_{m}^{u}-h_m+c_{mm}+c_{nm})^{2}}{4c_{mm}}\\
&\quad&\ds+(r_{n}^{u}-h_n+c_{nn})\varepsilon-\sum_{i=1}^{2}\int_{0}^{1}b_{i}(t)\mathrm{d}t\\
&\leqslant&\ds(r_{n}^{u}-h_n+c_{nn})\varepsilon+\frac{(r_{m}^{u}-h_m+c_{mm}+c_{nm})^{2}}{4c_{mm}}-\sum_{i=1}^{2}\int_{0}^{1}b_{i}(t)\mathrm{d}t\\
&\leqslant&\ds\frac{1}{2}\Big\{\frac{(r_{m}^{u}-h_m+c_{mm}+c_{nm})^{2}}{4c_{mm}}-\sum_{i=1}^{2}\int_{0}^{1}b_{i}(t)\mathrm{d}t\Big\}\\
&=&\ds-\frac{1}{2}\Phi_m\\
&\leqslant&\ds-1,
\end{array}
\end{equation*}
where $m=1,2,\quad m\neq n$.

Case 2: For any fixed $n (n=1,2)$, if $\frac{1}{\varepsilon}<x_{n}$, we have
\begin{equation*}
\begin{array}{rcl}
\mathcal{L}V&\leqslant&\ds-\frac{c_{nn}}{2}x_{n}^{2}(t)-\frac{c_{nn}}{2}x_{n}^{2}(t)-c_{mm}x_{m}^{2}(t)+\sum_{i=1}^{2}(r_{i}^{u}-h_i+c_{ii})x_{i}(t)\\
&\quad&\ds+\sum_{i,j=1,j\neq i}^{2}c_{ij}x_{j}(t)-\sum_{i=1}^{2}\int_{0}^{1}b_{i}(t)\mathrm{d}t\\
&\leqslant&\ds-\frac{c_{nn}}{2\varepsilon^{2}}+K\\
&\leqslant&-1,
\end{array}
\end{equation*}
where $m=1,2$, $m\neq n$ and $K=\sup_{(x_{1}, x_{2})\in R_+^2}\big\{-\frac{c_{nn}}{2}x_{n}^{2}(t)-c_{mm}x_{m}^{2}(t)+\sum_{i=1}^{2}(r_{i}^{u}-h_i+c_{ii})x_{i}(t)+\sum_{i,j=1,j\neq i}^{2}c_{ij}x_{j}(t)-\sum_{i=1}^{2}\int_{0}^{1}b_{i}(t)\mathrm{d}t\big\}$ is a constant.

In summary, we get
\begin{equation}\label{lv2}
\mathcal{L}V\leqslant-1,\quad(x_{1}, x_{2})\in R_+^2\setminus U_{\varepsilon} .
\end{equation}
From (\ref{lv1}) and (\ref{lv2}), we can note that the conditions in Theorem 3.8 of \citep{MR13} are satisfied. Therefore, model (\ref{fxx}) possesses a $1$-periodic solution which is denoted as $x^{*}(t)=(x_{1}^{*}(t), x_{2}^{*}(t))\in R_+^2$.
\end{proof}
\medskip

\rm
\section{Stability in distridution}
\label{4}

In this section, we consider the stability in distribution of model (\ref{fxx}). Before we prove our results, we state an assumption and a lemma.
\medskip

\noindent$\mathbf{Assumption\ 3.}$ $c_{11}>c_{21}$ and $c_{22}>c_{12}$.
\medskip

\begin{lemma}
If $p\geqslant1$, there exists a positive $K_2$ such that
\begin{equation*}
\limsup_{t\rightarrow+\infty}\mathbb{E}(x_i(t))^p\leqslant K_2, \quad i=1,2.
\end{equation*}
\end{lemma}
\medskip

\begin{proof}[$\mathbf{Proof\ of\ Lemma\ 4.1}$] \rm The proof of this lemma is a special case of lemma 3.1 in \citep{MR08}, so the proof process is omitted.
\end{proof}
\medskip

\begin{lemma}
If Assumption 3 holds, for any initial data $(x_1(0),x_2(0))^{T}$ and $(\widetilde{x}_1(0),\widetilde{x}_2(0))^{T}$, the solutions $(x_1(t),x_2(t))^{T}$ and $(\widetilde{x}_1(t),\widetilde{x}_2(t))^{T}$ obey
\begin{equation*}
\lim_{t\rightarrow+\infty}\mathbb{E}\big|x_i(t)-\widetilde{x}_i(t)\big|=0, \quad i=1,2.
\end{equation*}
\end{lemma}
\medskip

\begin{proof}[$\mathbf{Proof\ of\ Lemma\ 4.2}$] \rm Set
\begin{equation*}
\begin{array}{rcl}
\ds \widetilde{V}(t)\overset{\triangle}{=}\ds\big|\ln x_1(t)-\ln \widetilde{x}_1(t))\big|+\big|\ln x_2(t))-\ln \widetilde{x}_2(t))\big|.\\
\end{array}
\end{equation*}
Applying It$\rm\hat{o}$'s formula to compute the right differential $\mathrm{d}^+\widetilde{V}(t)$ of $\widetilde{V}(t)$, we obtain
\begin{equation*}
\begin{array}{rcl}
&\quad&\ds \mathrm{d}^+\widetilde{V}(t)\\&=&\ds \sgn\big(x_1(t)-\widetilde{x}_1(t))\big)\big[-c_{11}\big(x_1(t)-\widetilde{x}_1(t))\big)-c_{12}\big(x_2(t)-\widetilde{x}_2(t))\big)\big]\mathrm{d}t\\
&\quad&+\ds \sgn\big(x_2(t)-\widetilde{x}_2(t))\big)\big[-c_{21}\big(x_1(t)-\widetilde{x}_1(t))\big)-c_{22}\big(x_2(t)-\widetilde{x}_2(t))\big)\big]\mathrm{d}t\\
&\leqslant&\ds -\sum_{i=1}^{2}c_{ii}\big|x_i(t)-\widetilde{x}_i(t))\big|\mathrm{d}t\\
&\quad&+\ds
c_{12}\big|x_2(t)-\widetilde{x}_2(t))\big|\mathrm{d}t+c_{21}\big|x_1(t)-\widetilde{x}_1(t))\big|\mathrm{d}t.\\
\end{array}
\end{equation*}
Therefore,
\begin{equation*}
\begin{array}{rl}
\ds0\leqslant\mathbb{E}(\widetilde{V}(t))\leqslant&\ds \widetilde{V}(0)-(c_{11}-c_{21})\int_{0}^{t}\mathbb{E}\big|x_1(t)-\widetilde{x}_1(t))\big|\mathrm{d}s\\
\quad&\ds-(c_{22}-c_{12})\int_{0}^{t}\mathbb{E}\big|x_2(t)-\widetilde{x}_2(t))\big|\mathrm{d}s,
\end{array}
\end{equation*}
which implies that
\begin{equation*}
\begin{array}{rcl}
\ds(c_{11}-c_{21})\int_{0}^{t}\mathbb{E}\big|x_1(t)-\widetilde{x}_1(t))\big|\mathrm{d}s\leqslant \widetilde{V}(0)<+\infty,\\
\\
\ds(c_{22}-c_{12})\int_{0}^{t}\mathbb{E}\big|x_2(t)-\widetilde{x}_2(t))\big|\mathrm{d}s\leqslant \widetilde{V}(0)<+\infty.\\
\end{array}
\end{equation*}
Due to Assumption 3, we get, for every $t\in(0,\infty)$,
\begin{equation*}
\begin{array}{rcl}
\ds\int_{0}^{t}\mathbb{E}\big|x_1(t)-\widetilde{x}_1(t))\big|\mathrm{d}s\leqslant \frac{\widetilde{V}(0)}{(c_{11}-c_{21})}<+\infty,\\
\\
\ds\int_{0}^{t}\mathbb{E}\big|x_2(t)-\widetilde{x}_2(t))\big|\mathrm{d}s\leqslant \frac{\widetilde{V}(0)}{(c_{22}-c_{12})}<+\infty.\\
\end{array}
\end{equation*}
Thus,
\begin{equation}\label{w}
\mathbb{E}\big|x_i(t)-\widetilde{x}_i(t))\big|\in L^1[0,+\infty),\quad i=1,2.
\end{equation}
From model (\ref{fxx}), we note that
\begin{equation*}
\begin{array}{rcl}
\ds\mathbb{E}(x_1(t))&=&\ds x_1(0)+\int_{0}^{t}\Big[\mathbb{E}((r_1(t)-h_1)x_1(s))-c_{11}\mathbb{E}(x_1(s))^2\\
&\quad&\ds-c_{12}\mathbb{E}(x_1(s)x_2(s))\Big]\mathrm{d}s,
\end{array}
\end{equation*}
which implies the differentiability of $\mathbb{E}(x_1(t))$. Due to Lemma 4.1,
\begin{equation*}
\begin{array}{ll}
\ds \frac{\mathrm{d}\mathbb{E}(x_1(t))}{\mathrm{d}t}&\leqslant\ds (r_1^u-h_1)\mathbb{E}(x_1(t))-c_{11}\mathbb{E}(x_1(t))^2-c_{12}\mathbb{E}(x_1(t)x_2(t))\\
&\leqslant\ds r_1^u\mathbb{E}(x_1(t))\\
&\leqslant\ds r_1^uD_{1},
\end{array}
\end{equation*}
where $D_{1}>0$ is a constant. Thus $\mathbb{E}(x_1(t))$ is uniformly continuous. In a same way, $\mathbb{E}(x_2(t))$ is also uniformly continuous. By Barbalat's lemma (Barbalat \citep{MR14}) and (\ref{w}), we conclude:
\begin{equation}\label{e1}
\lim_{t\rightarrow+\infty}\mathbb{E}\big|x_i(t)-\widetilde{x}_i(t)\big|=0, \quad i=1,2.
\end{equation}
\end{proof}
\medskip

\begin{lemma}
If Assumption 3 holds, there is a unique $1$-periodic probability measure $\varphi(t, \cdot)$, such that, for every initial data $x(0)\in R^{2}_{+}$, the transition probability $p(t,0,x(0),\cdot)$ of $x(t)$ converges weakly to $\varphi(t, \cdot)$ as $t\rightarrow+\infty$.
In addition,
\begin{equation*}
\lim_{t\rightarrow+\infty}\frac{1}{t}\int^{t}_{0}x(s)\mathrm{d}s=\int^{1}_{0}\int^{+\infty}_{0}x\varphi(t, \mathrm{d}x)\mathrm{d}t, \quad a.s.
\end{equation*}
\end{lemma}
\medskip

\begin{proof}[$\mathbf{Proof\ of\ Lemma\ 4.3}$] \rm By Lemma 4.1, the $p$-th moment of $x_i(t)$ is uniformly bounded:
\[
\sup_{t \geq 0} \mathbb{E}[x_i(t)^p] \leq K \quad (i = 1, 2),
\]
where $K > 0$ is constant. Combined with inequality (23) and Theorem 2.8 in \citep{MR15}, the system admits a unique 1-periodic probability measure $\varphi(t, \cdot)$ such that for all $x(0) \in \mathbb{R}_+^2$, $p(t, 0, x(0), \cdot)$ converges weakly to $\varphi(t, \cdot)$.

Define the averaged measure $\overline{\varphi}$ on $\mathbb{R}_+$:
\[
\overline{\varphi}(\cdot) \overset{\triangle}{=} \int_0^1 \varphi(t, \cdot) \, \mathrm{d}t.
\]
By Theorem 3.2 in Feng and Zhao \citep{MR16}, $\overline{\varphi}(\cdot)$ is an invariant measure.

For $t \geq 0$, let $n_t$ be the maximal nonnegative
integer less or equal to $t$ and decompose:
\[
\frac{1}{t} \int_0^t x(s) \, \mathrm{d}s = \underbrace{\frac{n_t}{t} \cdot \frac{1}{n_t} \sum_{k=0}^{n_t-1} \int_k^{k+1} x(s) \, \mathrm{d}s}_{(I)} + \underbrace{\frac{1}{t} \int_{n_t}^t x(s) \, \mathrm{d}s}_{(II)}.
\]

\noindent \textbf{Term (I):} By the Markov property and periodicity, $\int_k^{k+1} x(s) \, \mathrm{d}s \overset{\triangle}{=} H(\theta_k \omega)$, where $\theta_k$ is the shift operator and $H(\omega) \overset{\triangle}{=} \int_0^1 x(s, \omega) \, \mathrm{d}s$.
Due to the uniqueness of $\varphi(t, \cdot)$ and the definition of $\overline{\varphi}$, as well as the properties of invariant measures, it can be concluded that $\bar{\varphi}$ is the unique invariant measure for the discrete system $\{H(\theta_k \omega)\}$. Consequently, $\bar{\varphi}$ is ergodic. The Birkhoff Ergodic Theorem 
yields:
\[
\lim_{n_t \to \infty} \frac{1}{n_t} \sum_{k=0}^{n_t-1} H(\theta_k \omega) = \mathbb{E}_{\overline{\varphi}}[H] \quad \text{a.s.}
\]

\noindent \textbf{Term (II):} By Lemma 4.1's moment bound:
\[
\left| \frac{1}{t} \int_{n_t}^t x(s) \, \mathrm{d}s \right| \leq \frac{t - n_t}{t} \cdot \sup_{s \geq 0} \mathbb{E}[x(s)] \to 0 \quad \text{as} \quad t \to \infty.
\]
By Fubini's theorem and $\overline{\varphi}$'s definition:
\begin{equation*}
\begin{array}{rl}
\mathbb{E}_{\overline{\varphi}}[H]=& \ds\int_{\Omega} H(\omega) \, \overline{\varphi}(\mathrm{d}\omega) \\
=& \ds\int_{\Omega} \left( \int_{0}^{1} x(s,\omega) \, \mathrm{d}s \right) \overline{\varphi}(\mathrm{d}\omega) \\
=& \ds\int_{0}^{1} \left( \int_{\Omega} x(s,\omega) \, \overline{\varphi}(\mathrm{d}\omega) \right) \mathrm{d}s \\
=& \ds\int_{0}^{1}\int_{0}^{+\infty} x \, \varphi(s,\mathrm{d}x)\mathrm{d}s.
\end{array}
\end{equation*}
Due to the definition of \( n_t\), \(n_t/t \to 1\) as \( t \to +\infty \).
Thus,
\begin{equation*}
\lim_{t\rightarrow+\infty}\frac{1}{t}\int^{t}_{0}x(s)\mathrm{d}s=\int^{1}_{0}\int^{+\infty}_{0}x\varphi(t, \mathrm{d}x)\mathrm{d}t, \quad a.s.
\end{equation*}
\end{proof}
\medskip

\rm
\section{Optimal harvesting}
\label{5}

In this section, we will state the optimal harvesting effort (OHE) and maximum sustainable yield (MESY) for model (\ref{fxx}).
\medskip

\begin{remark} Define
\begin{equation*}\quad
C=\left(\begin{array}{l}
 c_{11} \quad c_{12}  \\
 c_{21} \quad c_{22}  \\
\end{array}
\right).
\end{equation*}
It is necessary to point out that $C^{-1}+(C^{-1})^{T}$ is a positive definite matrix due to $c_{22}>0$ and $\Delta>0$.
\end{remark}
\medskip

\begin{theorem}
For model (\ref{fxx}), suppose that Assumptions 1-3 hold. Define
$$ L\overset{\triangle}{=}\Big(\int^{1}_{0}(r_1(t)-\frac{\alpha_{1}^{2}(t)}{2})\mathrm{d}t,\int^{1}_{0}(r_2(t)-\frac{\alpha_{2}^{2}(t)}{2})\mathrm{d}t\Big)^T,$$
\begin{equation}\label{abc}
\ds A\overset{\triangle}{=}(\lambda_1,\lambda_2)^{T}\overset{\triangle}{=}[C(C^{-1})^T+I]^{-1}L,
\end{equation}
where $I$ is a $2\times2$ identity matrix.

If $\int^{1}_{0}b_1(t)\mathrm{d}t|_{H=A}>0$, $\int^{1}_{0}b_2(t)\mathrm{d}t|_{H=A}>0$, $\Delta_1|_{H=A}>0$, $\Delta_2|_{H=A}>0$, $\Phi_1|_{H=A}>2$, $\Phi_2|_{H=A}>2$, $\lambda_1\geqslant0$ and $\lambda_2\geqslant0$.
Consequently, OHE is $H^*=A$ and MESY is
\begin{equation}\label{abd}
Y^*=A^{T}C^{-1}(L-A).
\end{equation}
\end{theorem}
\medskip

\begin{proof}[$\mathbf{Proof\ of\ Theorem\ 5.1}$]
\rm According to lemmas 3.1, 4.2 and 4.3, we have
\begin{equation}\label{xa}
\begin{array}{rl}
Y(H)=&\ds\liminf\limits_{t\rightarrow+\infty}\int_{t}^{t+1}\sum^2_{i=1}\mathbb{E}(h_ix_i(s))\mathrm{d}s\\
=&\ds\liminf\limits_{t\rightarrow+\infty}\int_{t}^{t+1}\mathbb{E}(H^Tx(s))\mathrm{d}s\\
=&\ds\liminf\limits_{t\rightarrow+\infty}\int_{t}^{t+1}H^T\mathbb{E}(x(s))\mathrm{d}s\\
=&\ds\liminf\limits_{t\rightarrow+\infty}\int_{t}^{t+1}H^T\mathbb{E}(x^*(s))\mathrm{d}s\\
=&\ds\int_{0}^{1}H^T\mathbb{E}(x^*(s))\mathrm{d}s\\
=&\ds\int^{1}_{0}H^T\int^{+\infty}_{0}x\varphi(t, dx)\mathrm{d}t\\
=&\ds\lim_{t\rightarrow+\infty}\frac{1}{t}\int^{t}_{0}H^Tx(s)\mathrm{d}s\\
=&\ds\sum^2_{i=1}h_i\lim_{t\rightarrow+\infty}\frac{1}{t}\int^{t}_{0}x_i(s)\mathrm{d}s.\\
\end{array}
\end{equation}
Combining (\ref{slu}) with (\ref{xa}), we obtian
\begin{equation}\label{xb}
\begin{array}{rl}
Y(H)
=&\ds \sum^2_{i=1}\frac{h_{i}\Delta_{i}}{\Delta}=\ds H^TC^{-1}(L-H).
\end{array}
\end{equation}
Let $A=(\lambda_1,\lambda_2)^{T}$ be the unique solution of the following equation:
\begin{equation*}
\begin{array}{rcl}
\ds 0=\frac{\mathrm{d}Y(H)}{\mathrm{d}H}&=&\ds\frac{\mathrm{d}(H^TC^{-1}(L-H))}{\mathrm{d}H}=\frac{\mathrm{d}(H^TC^{-1}L-H^TC^{-1}H)}{\mathrm{d}H}\\
&=&\ds\frac{\mathrm{d}(H^TC^{-1}L)}{\mathrm{d}H}-\frac{\mathrm{d}(H^TC^{-1}H)}{\mathrm{d}H}\\
&=&\ds C^{-1}L-[C^{-1}+(C^{-1})^T]H.
\end{array}
\end{equation*}
Hence $A=[C(C^{-1})^T+I]^{-1}L$. Clearly,
\begin{equation*}
\begin{array}{rcl}
\ds\frac{\mathrm{d}}{\mathrm{d}H^T}\Big[\frac{\mathrm{d}Y(H)}{\mathrm{d}H}\Big]&=&\ds\Big(\frac{\mathrm{d}}{\mathrm{d}H}\Big[(\frac{\mathrm{d}Y(H)}{\mathrm{d}H})^T\Big]\Big)^T\\
&=&\ds\Big(\frac{\mathrm{d}}{\mathrm{d}H}\Big[L^T(C^{-1})^T-H^T[C^{-1}+(C^{-1})^T]\Big]\Big)^T\\
&=&\ds-(C^{-1}+(C^{-1})^T).
\end{array}
\end{equation*}
is negative definite (see Remark 5.1)\color{black}. Therefore, $A=[C(C^{-1})^T+I]^{-1}L$ is the unique extreme point of $Y(H)$. Thus OHE $H^*=A$ and MESY is $A^{T}C^{-1}(L-A)$ by (\ref{xb}).
\end{proof}
\medskip

\rm
\section{Numerical simulations}
\label{6}

In this section, we will illustrate our main results by using the Euler-Maruyama method (seeing \citep{MR17}), leveraging Python for implementation. We always choose $r_1(t)=6.5+0.1\sin(2\pi t)$, $r_2(t)=6.6+0.1\sin(2\pi t)$, $c_{11}=4.3$, $c_{12}=0.4$, $c_{21}=0.5$, $c_{22}=3.5$ in this section. Then $\Delta=14.85>0$ and Assumptions 1, 3 hold.
\medskip

Firstly, we illustrate the effect of white noises on the optimal harvesting policy. We plot the curve of \( Y(H) \) in Fig. 1 by varying the parameters \( \alpha_1 \) and \( \alpha_2 \), under the initial conditions \( x_1(0) = 0.01 \) and \( x_2(0) = 0.01 \):

$\rm(\romannumeral1)$ The blue line is with $\alpha_1=0.1+0.01\cos(2\pi t)$, $\alpha_2=0.1+0.01\cos(2\pi t)$, then
$\int_{0}^{1}b_{1}(t)dt=3.20>0$, $\int_{0}^{1}b_{2}(t)dt=3.33>0$, $\Delta_1=9.88>0$, $\Delta_2=12.72>0$
$\Phi_1=2.71>2$,
$\Phi_2=2.69>2$, which means that
both $x_1$ and $x_2$ are persistence and Assumption 2 holds.
Hence, by Theorem 5.1, $H^*=(3.29,3.26)^T$, $Y^*=4.99$.

$\rm(\romannumeral2)$ The orange line is with $\alpha_1=0.7+0.01\cos(2\pi t)$, $\alpha_2=0.1+0.01\cos(2\pi t)$, then
$\int_{0}^{1}b_{1}(t)dt=3.08>0$, $\int_{0}^{1}b_{2}(t)dt=3.33>0$, $\Delta_1=9.46>0$, $\Delta_2=12.77>0$
$\Phi_1=2.48>2$,
$\Phi_2=2.57>2$, which means that
both $x_1$ and $x_2$ are persistence and Assumption 2 holds.
Hence, $H^*=(3.17,3.27)^T$, $Y^*=4.83$.

$\rm(\romannumeral3)$ The green line is with $\alpha_1=0.1+0.01\cos(2\pi t)$, $\alpha_2=1.1+0.01\cos(2\pi t)$, then
$\int_{0}^{1}b_{1}(t)dt=3.21>0$, $\int_{0}^{1}b_{2}(t)dt=3.03>0$, $\Delta_1=10.02>0$, $\Delta_2=11.43>0$
$\Phi_1=2.41>2$,
$\Phi_2=2.08>2$, which means that
both $x_1$ and $x_2$ are persistence and Assumption 2 holds.
Hence, $H^*=(3.29,2.96)^T$, $Y^*=4.50$.
%

\begin{figure}[htbp]
  \centering
  \includegraphics[width=1\textwidth]{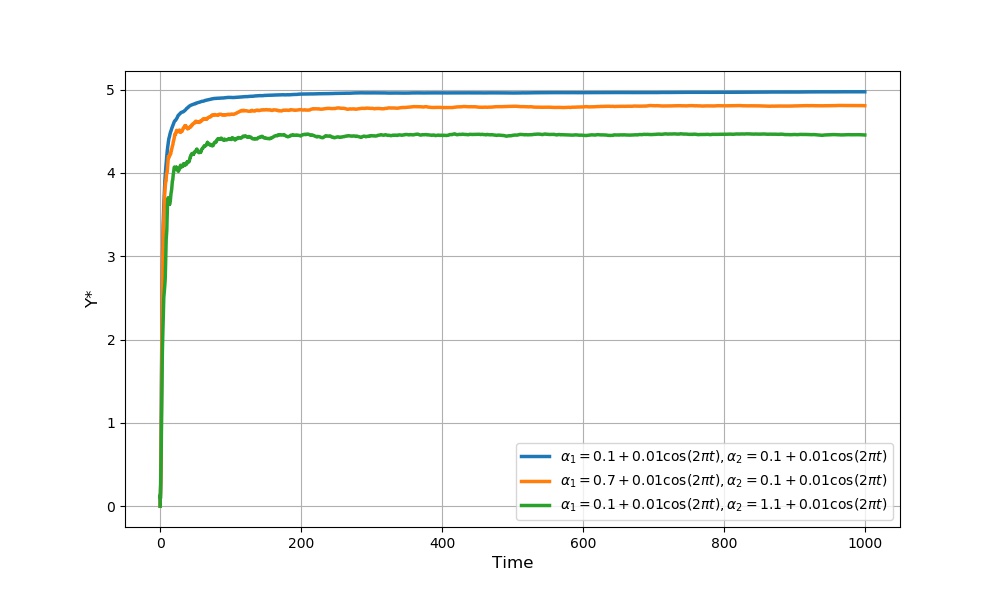}
  \small \footnotesize $\mathbf{Fig.\ 1.}$ The effect of white noises on $Y^*$. 
\end{figure}
\medskip

Finally, we illustrate that $H^*=[C(C^{-1})^T+I]^{-1}L$ and $Y^*=A^{T}C^{-1}(L-A)$ are the optimal harvesting policy in Fig. 2. We choose $\alpha_1=0.1+0.01\cos(2\pi t)$, $\alpha_2=0.1+0.01\cos(2\pi t)$, then $H^*=(3.29,3.26)^T$, $Y^*=4.99$.

\begin{figure}[htbp]
  \centering
  \includegraphics[width=1.1\textwidth]{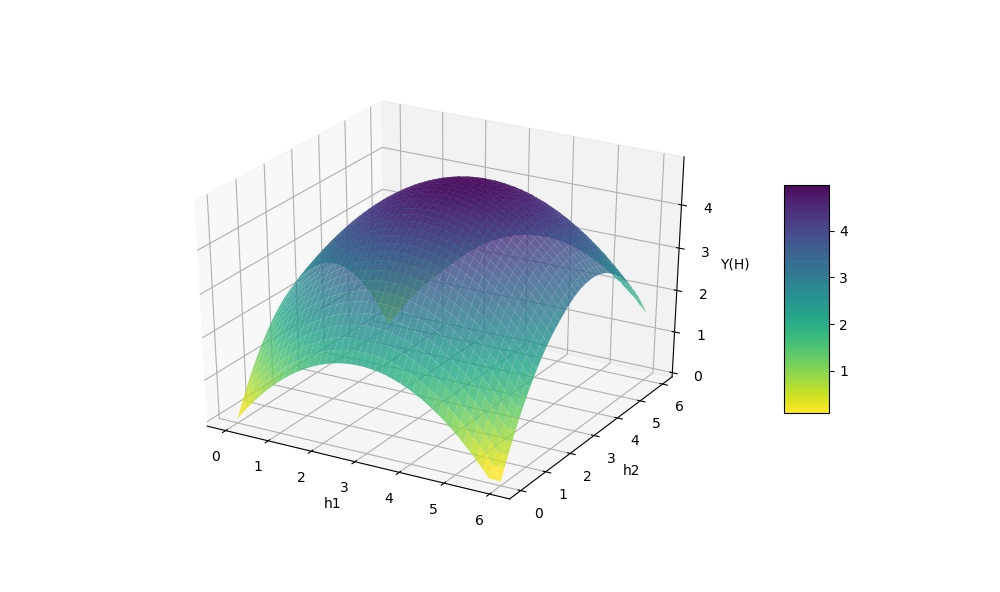}
  \small \footnotesize$\mathbf{Fig.\ 2.}$ Optimal Harvesting Policy $Y(H)$ as a function of $h_1$ and $h_2$. 
\end{figure}

\color{black}
\medskip

\rm
\section{Conclusions}
\label{7}

This paper investigates the optimal harvesting problem for a stochastic competitive Lotka-Volterra model with periodic coefficients. Within a rigorous mathematical framework, we establish necessary and sufficient conditions for the existence of an optimal harvesting policy (Theorem 5.1) and derive precise analytical expressions for OHE and MESY. Our results demonstrate a well-defined negative correlation between white noise intensity and optimal harvest quantities ((\ref{abc})-(\ref{abd})).
$$\frac{d\lambda_i}{d(\alpha^2_j(t))}\leqslant0, \quad \frac{dY^*}{d(\alpha^2_j(t))}\leqslant0,\quad i,j=1,2.$$

This finding carries significant ecological implications: environmental stochastic disturbances substantially increase species extinction risk, leading to decreased population abundance and consequently reducing both OHE and MESY. These results provide a theoretical foundation for biological resource management under stochastic environments.
\medskip

Based on this study, future research should focus on optimal harvesting strategies for stochastic population models with periodically varying coefficients in both interspecific and intraspecific interactions. It should be noted that since explicit expressions for time-averaged persistence conditions cannot be established for such models, the methodology presented in this study may not be directly applicable, thus requiring methodological improvements or alternative approaches.
\medskip

\rm
\label{9}


\end{document}